\documentclass{article}
\usepackage{amssymb,graphicx}

\def\d{{\rm d}}
\def\mi{{\rm i}}
\def\eps{\varepsilon}
\def\g{\gamma}
\def\G{\Gamma}
\def\l{\lambda}
\def\L{\Lambda}
\def\z{\zeta}
\def\Re{\mathop{\rm Re\,}\nolimits}
\def\Im{\mathop{\rm Im\,}\nolimits}
\def\e{\mathop{\rm e}\nolimits}

\def\hf{{\textstyle{1 \over 2}}}

\def\tq{{\textstyle{3 \over 4}}}

\def\defi{\stackrel{\rm def}{=}}
\def\si{\!\!\! &}
\def\se{& \!\!\!}
\def\ni{\noindent}

\newcommand{\beq}{\begin{equation}}
\newcommand{\eeq}{\end{equation}}
\newcommand{\bea}{\begin{eqnarray}}
\newcommand{\eea}{\end{eqnarray}}

\title{From asymptotic to closed forms\\
for the Keiper/Li approach\\
to the Riemann Hypothesis}
\author{{\bf Andr\'e Voros}\\
Universit\'e Paris--Saclay, CNRS, CEA, Institut de Physique Th\'eorique\\
91191 Gif-sur-Yvette, France\\
E-mail: {\tt andre.voros@ipht.fr}\\
dedicated to Yoshitsugu TAKEI for his $60^{\rm th}$ birthday}

\begin{document}

\maketitle

\begin{abstract}

The Riemann Hypothesis (RH) - that all nonreal zeros of Riemann's zeta function 
shall have real part 1/2 - remains a major open problem.
Its most concrete equivalent is that an infinite sequence of real numbers,
the Keiper--Li constants, shall be everywhere positive (Li's criterion).
But those numbers are analytically elusive and strenuous to compute, hence we seek simpler variants.
The essential sensitivity to RH of that sequence lies in its asymptotic tail;
then, retaining this feature, we can modify the Keiper--Li scheme
to obtain a new sequence in elementary closed form.
This makes for a more explicit analysis, with easier and faster computations.
We can moreover show how the new sequence will signal RH-violating zeros if any,
by observing its analogs for the Davenport--Heilbronn counterexamples to RH.

\end{abstract}

It is a great honor and pleasure to dedicate this talk to Professor Yoshitsugu Takei, 
for his major contributions to exact asymptotic analysis throughout his career since the early 90's, \cite{KT}
and surely for many more years to come. I am most grateful to the RIMS\footnote
{This workshop was supported by the Research Institute for Mathematical Sciences,
an International Joint Usage/Research Center located in Kyoto University.}
and the Organizers for their invitation - and many past ones.

After a digression on why and how (precisely 40 years ago!) we met exact asymptotic analysis, 
which was to be the source of our durable link and friendship with Y.~Takei (\S~1), 
we will mainly survey our recent work on the Keiper--Li approach to the Riemann Hypothesis 
- referring to \cite{V2} for any further detail. 
In \S~2 we review the original but elusive Keiper--Li sequence,
and then (\S~3) a discretization step (from derivatives to finite differences)
which leads to a modified sequence \emph{in elementary closed form}; 
its $n \to \infty$ behavior provides a new, very concrete, asymptotic criterion for the Riemann Hypothesis.
Finally, in \S~4 we test that criterion in generalized form, transposed to the Davenport--Heilbronn
functions which are counterexamples to the Riemann Hypothesis
(with some new material here: more data and discussions).
\bigskip

\emph{Riemann's zeta function} has the equivalent definitions, for $\Re x>1$: \cite{Ti}
\bea
\label{zd}
\z (x) &\defi& \sum_{k=1}^\infty k^{-x}
= \frac{1}{\G (x)} \int_0^\infty \frac{1}{\e^u-1} \, u^{x-1} \,\d u \\
\label{epf}
&=& \prod_{\{p\}} (1-p^{-x})^{-1} 
\mbox{ over all the prime numbers } p ;
\eea
the last form (Euler product) is just a quote to recall why $\z $ 
is a crucial function in number theory ($\log \z$ encodes the primes).

Standard actions upon the integral (Mellin) representation in (\ref{zd}) yield: 

\ni - that $\z$ is meromorphic in all of $\mathbb C$, 
with the only pole $\z (x) = \frac{\textstyle 1}{\textstyle x-1 \strut} + \cdots$\,;
 
\ni - the explicit values
\beq
\label{zbn}
\z (-n) = (-1)^n B_{n+1}/(n+1), \ n=0,1,2,\ldots \quad (B_n: \mbox{ Bernoulli numbers});
\eeq
- and Riemann's Functional Equation, best written as
\beq
\label{rfe}
2\xi(x) = 2\xi (1-x), \quad \mbox{where } 
\ 2\xi (x) \defi x(x-1) \pi^{-x/2}\G (x/2) \, \z (x) ;
\eeq
$2\xi$ is dubbed ``completed zeta function" ($\xi $ is the classic choice, 
but $2\xi$ is better normalized for us, with $2\xi (0)=2\xi (1)=1$).
The Functional Equation and (\ref{zbn}) imply the further explicit values
\beq
\label{z2m}
\z (2m) = \frac{|B_{2m}|}{2 (2m)!} \, (2\pi)^{2m} \iff 
2\xi (2m) = \frac{|B_{2m}|}{(2m \!-\! 3)!!} \, (2\pi)^m \quad (m=1,2,\ldots)
\eeq
where $k!! = 2^{(k+1)/2} \, \G (\hf k+1) / \sqrt \pi$
for $k$ odd (the usual double factorial).

\section{Aside: our first use of an exact WKB method}

Our strong link with Y. Takei and the Japanese school of complex analysis
stems from the growth of an exact form of asymptotics around 1980,
itself much inspired and encouraged by M. Sato, T. Kawai, M. Kashiwara. \cite{SKK}
(Other precursors included Leray, Boutet de Monvel--Kr\'ee, Bender--Wu, Dingle, Balian--Bloch, 
Sibuya, Zinn-Justin, as quoted in \cite[\S~1.2]{Ve}.) 
But since our main topic will only touch standard (as opposed to exact) asymptotics,
it may be timely here to share our personal recollection of why we met exact asymptotics at all, 
as this was by quite an accidental circumstance, not much told, 
and moreover sharing the preliminaries (\ref{zd})--(\ref{z2m}) of our later main topic: 
we wanted to generalize an asymptotic form of Riemann's Functional Equation (\ref{rfe}), 
and that needed an exact complex-WKB treatment!
\medskip

A pending problem in the 70's was the 1D quantum \emph{quartic} oscillator 
(the Schr\"odinger operator $-\d^2/\d q^2 + q^4$, $q \in \mathbb R$)
and specially its spectrum, only known to be discrete 
($\{E_\ell \}_{\ell=0,1,\ldots }$; $E_\ell >0$, $E_\ell \uparrow +\infty$) 
and to solve an eigenvalue condition of the form
${\mathcal S}(E_\ell) = 2\pi (\ell + \hf)$ for integer $\ell$
with a function ${\mathcal S}(E)$ supplied asymptotically:
\beq
\label{bsr}
{\mathcal S}(E) \sim b_0 E^{\frac{3}{4}} + \sum_{n=1}^\infty b_n E^{-\frac{3}{4}(2n-1)} 
\quad \mbox{for } E \to +\infty
\eeq
(Bohr--Sommerfeld), with $b_0 = \oint _{\{p^2+q^4=1\}} p \,\d q\ (>0) :$ a classical-action period.

In 1979, following \cite[\S~7]{MP}\cite{P}, we considered the \emph{spectral zeta function}
\beq
\label{szf}
Z(x) \defi \sum_\ell E_\ell ^{\, -x} \qquad (\Re x > \tq):
\eeq
a Dirichlet series like (\ref{zd}) for $\z (x)$, but with the \emph{thoroughly unknown}
eigenvalues $E_\ell $ in place of the integers $k$. Remarkably though, 
$Z(x)$ kept many (of the non-arithmetic) \emph{explicit} properties present in Riemann's $\z (x)$: 
\cite{VGS}\cite{Vz}

\ni - a Mellin representation, implying that $Z(x)$ is meromorphic in all of $\mathbb C$
and all its polar singularities can be written out;

\ni - explicit finite values: all $Z(-n), \ (n=0,1,2,\ldots)$, \cite{P} plus $Z'(0)$ and $Z(1)$.

There is just \emph{no functional equation} for $Z(x)$ to generalize Riemann's eq.~(\ref{rfe}) for $\z (x)$ 
(which links to the \emph{harmonic} ($q^2$) oscillator, of spectral zeta function $(1-2^{-x})\z (x)$). 
But if we only watch $x \to \pm \infty $ asymptotics, then ${\z (x) \sim 1}$ for $x \to +\infty $
reducing (\ref{rfe}) to: $\z (x) \sim \G (1-x)(2\pi)^x \sin(\hf \pi x)/\pi$ for $x \to -\infty$.
Now this remnant of (\ref{rfe}), just an \emph{explicit} $(x \to -\infty)$ asymptotic formula,
may generalize to other Mellin transforms: such a function, $\int \Theta (u) u^{x-1} \d u$,
\emph{potentially} has its $x \to -\infty$ behavior dictated, and thus described,
by the \emph{nearest singularities of} $\Theta (u)$ in ${\mathbb C}^\ast$
- i.e., \emph{provided the latter are isolated and computable}. 
And all that worked for $Z(x)$, under the specific Mellin representation

\begin{figure}[h]
\center 
\includegraphics[scale=.25]{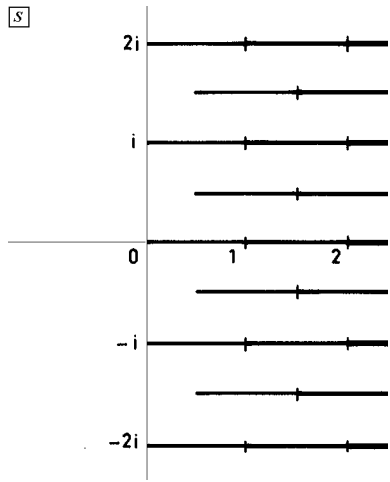}
\caption{\label{fg1} \small Cut plane (principal branch) for the multivalued function $\Theta_{3/4}(u)$, 
in the rescaled and negated variable $s=-b_0^{\,-1} u$. (Based on \cite[Fig.~11]{BPV}.)}
\end{figure}

\beq
~Z \bigl( \tq x \bigr) = 
\frac{1}{\G (x)} \int_0^\infty \Theta_{3/4}(u) \, u^{x-1} \d u, 
\quad \Theta_{3/4}(u) \defi \sum_\ell \exp(-E_\ell ^{\,3/4} u) \ \ (\Re u>0),
\eeq
because \emph{this} function $\Theta_{3/4}$ had a curious (and novel, at the time) analytic structure
depicted in \cite[\S~4]{BPV}, and here in Fig.~\ref{fg1}, using a rescaled variable $s$: 
$\Theta_{3/4}$ was a \emph{ramified} function, with branch points all on a square lattice, 
and up to the rescaling, its discontinuity functions were: 
at 0, a Borel transform of the Bohr--Sommerfeld series~(\ref{bsr}), 
and at other lattice points, Borel transforms of various 
\emph{exponentials of that same Bohr--Sommerfeld series} (essentially).
In particular that gave the nearest discontinuity functions, 
at $s=\hf (1 \pm \mi)$, to all orders, \cite[eq.(4.12)]{BPV} 
implying this asymptotic expansion for $Z(x)$, \cite[\S~V]{VGS}\cite[\S~7]{Vz}
\beq
\label{asz}
Z(\tq x) \sim \frac{\sin \tq \pi x}{\cos \hf \pi x} \biggl( \frac{b_0}{\sqrt 2} \biggr) ^{\! x} 
\frac{x}{\pi } \Biggl[ \sum_{j=0}^\infty \alpha _j \G (-x-j) \Biggr] \qquad (x \to -\infty),
\eeq
with $\{ \alpha _j\}$ ``bootstrapped" \emph{in terms of} $\{b_n\}$ through the generating function 
\beq
\label{bst}
\sum_{j=0}^\infty \alpha _j \Bigl(\frac{v}{b_0} \Bigr) ^j =
\exp \biggl[ \Bigl( \frac{b_1}{2} v + \frac{b_2}{2^2} v^3 \Bigr)
- \Bigl( \frac{b_3}{2^3} v^5 + \frac{b_4}{2^4} v^7 \Bigr) + \cdots \biggr] \quad (e.g., \ \alpha _0 = 1) .
\eeq

Yet we noted that our crucial description of the function $\Theta_{3/4}$
(later understood as \emph{resurgent} in the sense of \'Ecalle \cite{E}) 
stayed somewhat empirical and incomplete.
And however fast eq.~(\ref{asz}) grew for $x \to -\infty$, it still eluded a standard asymptotic approach. 
So, to confirm (\ref{asz})--(\ref{bst}) we needed tools able to \emph{fully} describe $\Theta_{3/4}$.
For that (encouraged also by Balian, Malgrange) we had to do WKB calculations 
with complex Planck's constant as in \cite{BB}, using microfunction techniques \cite{SKK}
and a whole convolution algebra of Borel transforms; thus, from 1981 onwards
\cite[\S~5.3]{Vc}\cite{V0}\cite{VB}\cite{V1} we ended up with \emph{exact WKB results}
- initially all for the sake of that $x \to -\infty$ behavior in the spectral zeta function $Z(x)$.
(For which \cite[pp.281--286]{V1} fully analyzed the function $\Theta_{3/4}$ - named $Z$ therein.)

As we also use (\ref{z2m}) later, we mention that even without a functional equation for $Z(x)$,
our further exact-WKB study of the potential $q^4$ extended the formulae 
(\ref{z2m}) for $\z (2m)$, to explicit identities for $Z(3m), \ m=1,2,\ldots$; 
e.g., $Z(3)=\frac{1}{6}Z(1)^3-\hf Z(1)Z(2)$. \cite{V0}
And similarly for higher-degree potentials $q^{2M}$ ($M>2$), and for parity-twisted zeta functions 
$\sum\limits_\ell (-1)^\ell E_\ell ^{\,-x}$ as well;
then, \emph{one} such zeta-value exceptionally reduces almost as far as~(\ref{z2m}):
\beq
\sum_\ell (-1)^\ell E_\ell ^{\,-2} = \textstyle \frac{1}{128} \, [\pi \, \G (\frac{1}{4})]^2/\G (\frac{7}{8})^4
\quad \mbox{for the \emph{sextic} potential } q^6
\eeq
by merging \cite[eq.~(16)]{V0} (exact-WKB) 
and \cite[eq.~(12)]{VGS}\cite{Vz} (Weber--Schafheitlin) at $\mu \defi (2M+2)^{-1} = 1/8$.

Still, in our own work, for us several gaps \emph{remain to be filled}: 
e.g., proof of full resurgence for the WKB solutions \cite{KS}, 
validation of a very gentle behavior at infinity seen on Borel transforms \cite[p.102--103]{VB},
regularity and contractivity of an exact-quantization map  
beyond the case of homogeneous potentials which was settled by Avila 
\cite{A} (\cite[\S~4.2]{Ve}\cite[end~\S]{VK}), 
broader inclusion of nonpolynomial potentials, \cite[\S~5.3]{Ve}
exact WKB treatment in phase-space quantum mechanics,\ldots

We conclude this digression by a salute to the impressive developments further carried out 
on exact WKB analysis in Japan, encompassing higher-order Ordinary Differential Equations, 
singular potentials, infinite-dimensional problems, ODE/IM correspondence,
nonlinear problems (Painlev\'e), cluster algebras\ldots

\section{The Keiper--Li sequence}

\subsection{The Riemann zeros (basics) \cite{Ti}}

\subsubsection{Known facts}

By (\ref{rfe}), $\xi $ is a real entire function,
with the two symmetry axes $\mathbb R$, and $L \defi {\{\Re x=\hf \}}$
called the \emph{critical line}.

The zeros of $\xi $ or \emph{Riemann zeros}, classically denoted $\rho $ 
(and counted with multiplicities if any), all lie within the open strip $\{0<\Re x<1\}$.
They are infinitely many and their counting function $N(T)$, defined as the number of $\rho$ 
in the rectangle $(0,1)\times(0,\mi T)$, obeys the \emph{Riemann--von Mangoldt} asymptotic law
\beq
\label{rvm}
N(T) = \frac{T}{2\pi} \Bigl( \log \frac{T}{2\pi}-1 \Bigr) +
O(\log T) \qquad (T\to +\infty).
\eeq

\subsubsection{The Riemann Hypothesis (RH) (1859) \cite{Ri}}

\centerline{\emph{All the zeros $\rho $ of $\xi(x)$ lie 
on the critical line} $\{ \Re x = \hf \}$.}
\smallskip

This conjecture, most important for number theory (to understand the primes) 
has been neither proved nor disproved yet. On the other hand:

\ni - $\Re \rho = \hf$ has been seen, and verified by computer, 
up to increasing ordinates $T=\Im \rho \,$: since 2004, up to the $10^{13}$-th zero $\rho \,$; \cite{G}
that sets the largest ordinate $T_0$ up to which RH is verified to a current value $\approx 2.4 \cdot 10^{12}$.

\ni - numerous statements equivalent to RH, or criteria for RH, have been issued; 
many are highly abstract, but our focus will be on a specially concrete and simple-looking one.

\subsection{The Keiper--Li tool to test RH}

\subsubsection{The Keiper vs Li sequences: generalities}

Those sequences are defined: by the generating function (Keiper \cite{K})
\bea
\label{lnk}
\sum_{n=1}^\infty \l _n^{\rm K} z^n \si=\se \Phi(z) \defi
\log 2\xi \Bigl( x = \frac{1}{1-z} \Bigr) \\
\label{lkn}
\iff \l _n^{\rm K} \si=\se
\frac{1}{2 \pi \mi} \oint_C \frac{\d z }{z^{n+1}} \, \Phi (z) ,
\ \ C = \{|z|=\eps \ll 1\} \mbox{ positively oriented}; \quad
\eea
resp. by sums over all Riemann zeros grouped symmetrically (Li \cite{LI1})
\bea
\label{lnl}
\l _n^{\rm L} &\defi& \sum_\rho [1-(1-1/\rho)^n], \qquad n=1,2,\ldots \\
\label{lcl} 
&=& \sum_\rho [1-\cos n \theta_\rho], \qquad
\theta _\rho \defi -\mi \log (1-1/\rho).
\eea
Both are denoted $\l _n$ in the literature, but beware: 
$\l _n^{\rm L} = n\l _n^{\rm K}$; so in way of a pun, 
Keiper's $\l _n$ and Li's $\l _n$ \emph{differ by their common notation}.
Neither normalization is nicer on all counts, so we rather keep both
and use disambiguation superscripts K, L when the factor $n$ matters.

\begin{figure}[h]
\center
\includegraphics[scale=.4]{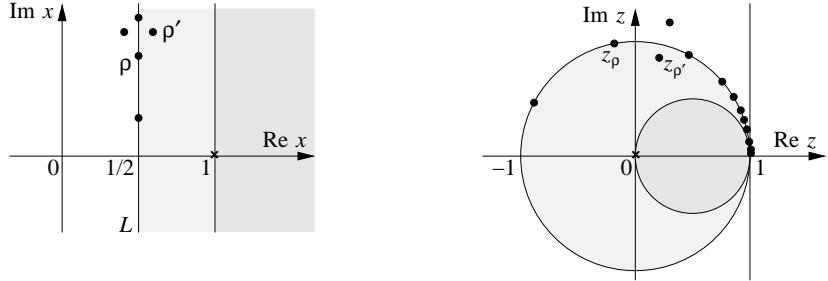}
\caption{\label{fg2} \small Riemann zeros ({\footnotesize \textbullet}) depicted 
in the $x$ (left) and $z$ (right) upper half-planes \emph{schematically} 
(at mock locations, including a putative pair off the critical line~$L$).
The symmetrical zeros in the lower half-planes are not plotted.
Domains are shaded only to mark which map to which.}
\end{figure} 

A key ingredient is the conformal mapping $x=1/(1-z)$ in (\ref{lnk}) 
which pulls back the half-plane $\{\Re x>\hf \}$ to the unit disk $\{|z|<1\}$ (Fig.~\ref{fg2}).
This makes RH equivalent to: $\Phi (z)$ is analytic in \emph{all of} that disk
- that is why its Taylor coefficients $\l _n$ are RH-sensitive. 
In quantitative terms, $\d \Phi /\d z$ is meromorphic with a simple pole of residue 1 
at every preimage $z_\rho = 1 - 1/\rho $ of a zero $\rho$, and by (\ref{lkn}),
\beq
\label{lll}
\l _n^{\rm L} = \frac{1}{2 \pi \mi} \oint_C \frac{\d z}{z^n} \,\frac {\d\Phi }{\d z} .
\eeq

\begin{figure}[h]
\center
\includegraphics[scale=.4]{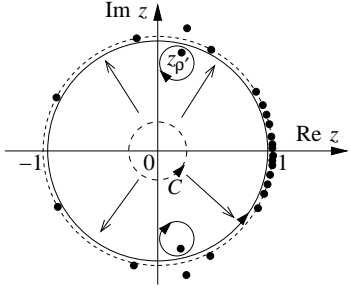}
\caption{\label{fg3} \small Contour deformation in Darboux's method for the asymptotics of (\ref{lll}).}
\end{figure}

\ni Applying Darboux's method (i.e., the method of steepest descent in the variable $\log z$ \cite[\S~7.2]{Di})
we inflate $C$ to $\{|z|=r \}$ with $r \uparrow 1^-$ (Fig.~\ref{fg3}) and use the residue theorem 
to get the contributions from the poles $z_{\rho '}$ as 
\beq
\label{rhf}
\l _n^{\rm L} = -\sum_{\{ |z_{\rho '}|<1\}} \! z_{\rho '}^{\, -n} 
\quad {}+ o(r^{-n})_{n \to \infty} \quad (\forall r<1)
\eeq
where we assign the notation $\rho '$ to zeros (if any) having $\Re \rho '> \hf$ 
(in violation of RH, and amounting to $|z_{\rho '}|<1$).
If and only if RH is false, the sum in (\ref{rhf}) is nonempty and then,
ordered according to nondecreasing $|z_{\rho '}|$ it forms an asymptotic expansion 
in \emph{exponentially growing} oscillations \emph{about $0$}.

\subsubsection{Li's criterion for the Riemann Hypothesis}

\ni - If RH is false, the last sentence about (\ref{rhf}) implies that 
$\l _n<0$ \emph{will} occur in the asymptotic regime $n \to \infty$.

\ni - If RH is true, this amounts to all $\theta_\rho $
being real in the sums (\ref{lcl}), which are therefore termwise \emph{positive}, for all $n$. \cite{K}

That pair of statements entails \emph{Li's criterion}: \cite{LI1}
\smallskip

\centerline{{\bf RH true} $\ \iff \ \l _n >0$ for all $n$.}

However: \cite{O}

\centerline{$\Re \rho = \hf$ holds up to a height $T_0$ \quad $\Longrightarrow \quad 
\l _n>0$ as long as $n<T_0^{\, 2}$.}

This means that low values of $n$ are actually inessential for Li's criterion: 
we may focus on the asymptotic $n \to \infty$ behavior of $\l _n$ instead.

\subsection{Asymptotic alternative for RH}
\label{aal}

The $n \to \infty$ form of $\l _n$ is already fixed by the sum in (\ref{rhf}) for RH false,  
but not so for RH true when that sum is empty.
Instead, in the RH true case the sum (\ref{lcl}) defining $\l _n^{\rm L}$ identifies with the
Stieltjes integral $2 \int_0^\infty (1-\cos n \theta) \, \d N(T)$ where $T = \hf \cot \hf \theta$;
then, integration by parts gives
\beq
\label{osf}
\l _n^{\rm K} = 2 \int_0^\pi \sin n \theta \, N(\hf \cot \hf \theta) \, \d \theta .
\eeq
The large-$T$ law (\ref{rvm}) now gives the $\theta \to 0$ form of the integrand in (\ref{osf}),
which in turn converts to the large-$n$ behavior of $\l _n^{\rm K}$ for RH true, as \cite{O}
\beq
\label{rht}
\l _n^{\rm K} = \hf [\log n + (\g - \log 2\pi -1)] +o(1) \qquad (\g : \mbox{ Euler's constant}),
\eeq
giving a \emph{tempered growth to} $+\infty$ (see \cite{Lg}\cite{AR} for stronger remainder estimates).
\smallskip

The asymptotic forms (\ref{rhf}), (\ref{rht}) together imply 
the ($n \to \infty$) \emph{alternative} \cite{V}
\bea
\label{rf}
\l _n^{\rm L} &\sim& 
- \! \sum\limits_{\{ |z_{\rho '}|<1\}} \! z_{\rho '}^{\, -n}
 \qquad \qquad \qquad \qquad \ \mbox{\bf if RH false} \\
\label{rt}
\mbox{\bf vs }\qquad \l _n^{\rm L} &\sim& \hf n \, [\log n + (\g - \log 2\pi -1)] \qquad \mbox{\bf if RH true.} 
\eea

In practice, a term $z_{\rho '}^{\, -n}$ from (\ref{rf}) will compete in size with (\ref{rt}) if
\beq
\label{upl}
n \gtrsim T^2/t \qquad (\mbox{for } \rho '= \hf + t \pm \mi T , \quad t>0).
\eeq
This inequality (in order of magnitude) is also the \emph{uncertainty principle} 
for the Fourier-conjugate variables $\theta $ and $n$ in (\ref{osf}),
which proves it a strict necessary condition as well. 
With $t<\hf$ and $T>T_0 \approx 2.4 \cdot 10^{12}$, (\ref{upl}) gives as concrete threshold: $n \gtrsim 10^{25}$,
for $\l _n$ to possibly sense violations of RH if any.

Unfortunately, the Keiper--Li numbers seem analytically quite challenging \cite{BL}\cite{C2}
and the complexity of their numerical evaluation steeply grows with~$n$,
mainly because it needs derivatives $(\log \xi )^{(n)}$ which are intricate to handle
\cite{K}\cite{M1}\cite{C1}\cite{J} (\cite{J} reached $n=10^5$).
Currently, only the behavior (\ref{rt}) will show over the accessed ranges (Fig.~\ref{fg4} below), 
whereas values $n \gtrsim 10^{25}$ needed for new tests of RH appear way out of reach.

\section{A closed-form variant of Keiper--Li \cite{V2}}

The construction of the Keiper--Li sequence is not as inflexible as it may seem.
Already, the unit disk of Fig.~\ref{fg2} can be remapped 
to itself by a (conformal) M\"obius transformation,
$z \mapsto z' = \frac{\textstyle z-\tilde z}{\textstyle 1-\tilde z^\ast z} \defi H_{\tilde z} (z)$
(for $|\tilde z|<1$; ${}^\ast$ denotes complex conjugation). 
Under the resulting composed map ${x \mapsto z \mapsto z'}$, 
$z'=0$ can now correspond to an arbitrary point $x_0$ in the half-plane $\{\Re x>\hf\}$, 
and the transposition of (\ref{lnk}) will define generalized Keiper--Li numbers 
in terms of ($\log \xi $)-derivatives now at $x=x_0$ instead of $x=1$
\cite{Se} (which improves convergence if $\Re x_0 >1$).
But $\log \xi $ stays differentiated all the same, only elsewhere.
To progress further, the deformation can be made more general.

\subsection{Construction of an explicit sequence}

The unwelcome differentiations on $\log \xi$ relate to the multiplicity
of the pole $1/z^{n+1}$ of the integrand in (\ref{lkn}).
So we propose to split this pole into $n+1$ simple poles 
$1/z,\ 1/(z-z_1), \ldots 1/(z-z_n)$
by displacing each factor of $z^{n+1}$ in (\ref{lkn}) differently.
As previously we use hyperbolic displacements, i.e., M\"obius transformations $H_{\tilde z}$, 
again to keep the unit disk invariant (thus preserving asymptotic sensitivity to RH). 
The integral form (\ref{lkn}) for \emph{Keiper's} $\l _n^{\rm K}$ thus gives rise to
\[
\frac{1}{2 \pi \mi} \oint_C \frac{\d z }{z \, H_{z_1}(z) \cdots H_{z_n}(z)} \, \Phi (z)
= \sum_{m=1}^n \frac{1}{z_m [H_{z_1}(z) \cdots H_{z_n}(z)]'(z_m)} \, \Phi _m
\]
by the elementary residue calculus for simple poles, with $\Phi _m \defi \Phi (z_m)$: 
i.e., a \emph{finite-difference} formula replaces a differential one.
Finally specializing to $z_m = 1-(2m)^{-1}$, the above reduces to
\bea
\label{ldd}
\L _n &\defi& \sum_{m=1}^n (-1)^m A_{nm} \Phi _m , \qquad \mbox{with} \\
A_{nm} &=& \frac{2^{m-n} \, \bigl( 2(n+m)-1 \bigr) !!}{(2m-1) \, (n-m)! \, (2m)!}, \\
\mbox{and} \quad \Phi _m &=& \log 2 \xi (2m) = 
\log \biggl[ \frac{|B_{2m}| }{(2m \!-\! 3)!! } \, (2\pi)^m \biggr]
\eea
thanks to (\ref{z2m}); thus, (\ref{ldd}) specifies a deformed Keiper's sequence $\{\L _n\}$ 
in elementary closed form. E.g.,
$\L _1 = \frac{3}{2} \log \pi/3 \approx 0.0691764, 
\ \L _2 = \frac{5}{24} \log [(2/5)^7 3^{11}/\pi^4 ] \approx 0.2274543$~.

In summary: the original $\l _n$ are elusive objects partly because 
the functions $\z $ hence $\xi $, $\log \xi $, have unwieldy derivatives. 
By discretizing the latter to finite differences, we inversely use the best 
in those same functions: their special values - countably many -
upon which explicit finite differences can be built - to all orders, just as needed here.

\subsection{Asymptotic alternative for RH with the sequence $\{\L _n\}$}

We only state our main result, in parallel to \S~\ref{aal}: defining
\beq
F_n (x) \defi
(-1)^n \biggl[ - \frac{1}{A_{n0}} \log (x-1) + \sum_{m=0}^n (-1)^m A_{nm} \log (x-2m) \biggr] 
\eeq
as single-valued in the cut plane ${\mathbb C} \setminus [0,2n]$,
then for $n \to \infty$,
\beq
\label{hf}
\qquad \qquad \ \L _n \sim \! \sum_{\{ \Re \rho '>1/2 \} } \! F_n(\rho ') \qquad \qquad \qquad \qquad \qquad \quad \mbox{\bf if RH false}
\eeq
where for \emph{each} $ \rho '=\hf + t + \mi T, \ t>0$ (with $\phi(\rho ')$: a known phase function),
\bea
\label{AS}
F_n(\rho ') &\sim& \e^{\mi \phi(\rho ')} \frac{(-1)^n (2n)^{t+\mi T}}{|T|^{2+t} \log n} \qquad 
\mbox{for } n \gg |T| \\
&& \mbox{(\emph{giving an oscillation of amplitude $O(n^t/\log n)$ about $0$});} \nonumber\\[10pt]
\label{ht}
\mbox{\bf vs }\qquad \L _n &\sim& \log n + \hf (\g - \log \pi -1) \pmod{o(1)} \qquad \mbox{\bf if RH true} \quad \\
&& \mbox{(\emph{tempered growth to} $+\infty$).} \nonumber
\eea
(The derivations are similar to those sketched above for (\ref{rf})--(\ref{rt})
but more elaborate; the explicit form (\ref{AS}) doesn't hold uniformly in $T = \Im \rho '$
hence cannot be substituted all at once into the full sum (\ref{hf}). \cite[\S~3]{V2})

\begin{figure}[t]
\center
\includegraphics[scale=.39]{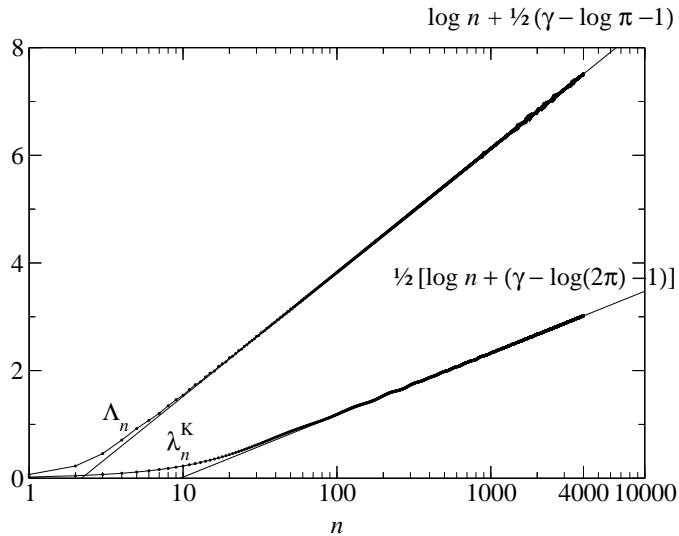}
\vskip -3mm

\caption{\label{fg4} \small Upper plot: the sequence $\{\L _n\}$ given by (\ref{ldd}) displayed up to $n=4000$
on a logarithmic $n$-scale (line segments connect data points only to aid the eye);
straight line: the RH-true asymptotic form (\ref{ht}). 
Lower plot: the same for the original Keiper sequence (\ref{lkn}) in comparison 
(data by courtesy of K. Ma\'slanka \cite{M1}).}
\end{figure}

\begin{figure}[h]
\center
\includegraphics[scale=.39]{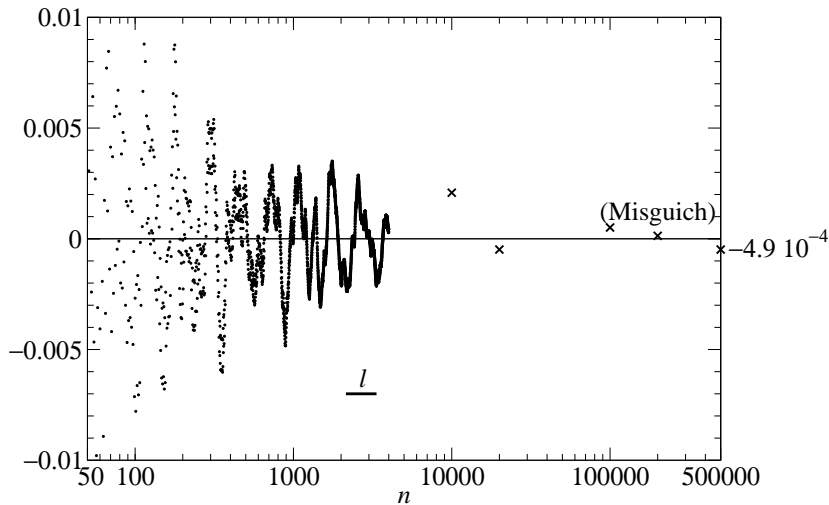}
\vskip -3mm

\caption{\label{fg5} \small $(-1)^n$ times the remainder $\delta \L _n$ 
in (\ref{ht}) ($=o(1)$ iff RH is true). 
This form oscillates with roughly the ($\log n$)-wavelength of $n^\rho $ 
for the first Riemann zero $\rho =\hf + 14.134725\ldots \, \mi$;
this wavelength, 
$l = 2\pi/\Im \rho \approx 0.44452$, is depicted by a horizontal segment.
Rightmost data: courtesy of G. Misguich (see \S~\ref{csp}).}
\end{figure}

In practice, a term from $\rho '$ as in (\ref{AS}) will compete in size with (\ref{ht}) if
\beq
\label{TN}
n \gtrsim T^{1+2/t} \qquad (\gtrsim T_0^{\,5} \approx 10^{60} \mbox{ currently): sufficient condition};
\eeq
now the uncertainty principle is much more favorable than either (\ref{TN}) or (\ref{upl}):
\beq
\label{ucp}
n \gtrsim \hf T \e ^{1/t} \qquad (\gtrsim \hf T_0 \e ^2 \approx 10^{13} \mbox{ currently): necessary condition}.
\eeq
Still, current data (Figs.~\ref{fg4}--\ref{fg5}) will then only show the behavior (\ref{ht}).

The two behaviors (\ref{hf}), (\ref{ht}) are mutually exclusive asymptotically (``alternative"), 
but numerically they superpose (they add): the form (\ref{ht}) sums 
the bulk effect of the zeros on the critical line, but the remainder therein,
$\delta \L _n = \L _n - \log n - {\hf (\g - \log \pi -1)}$ 
retains oscillations probably due to those zeros taken individually:
e.g., the main oscillation is clearly synchronous with $(-1)^n n^\rho $, 
the same form as in (\ref{AS}) but for the first Riemann zero 
($\rho = \hf + 14.134725 \ldots \, \mi$) (Fig.~\ref{fg5}).
Inversely, any growing oscillation (\ref{AS}) from an RH-violating zero $\rho '$
will rise on top, not in place, of the smooth trend (\ref{ht}),
see later counterexample to RH (Fig.~\ref{fg6}).

\subsection{Computational aspects}
\label{csp}

Calculations on $\L _n$ appear much simpler than for $\l _n$. A handful of command lines suffice in Mathematica (for instance) \cite{W}
to readily obtain values up to $n=20,000$. G. Misguich has written a much faster parallel
code for a 20-core machine, reaching $\L _{500,000} \approx 12.33812102688$ in about 22 days. \cite{Mi}

As a bonus, individual $\L _n$ can be accessed directly without a recursive build-up
from $n=1$ every time as with $\l _n$ (evaluations of the Bernoulli numbers $B_{2m}$ use recursion though, 
causing the major part of the workload).

An important fact is that the observed small $\L _n$ follow from huge cancellations 
between positive and negative terms in (\ref{ldd}). 
Here, this too can be described explicitly:
for $n \gg 1$ the Stirling formula shows $\max\limits_m |A_{nm}|$
to grow like $\e ^{(3+2\sqrt 2)n}$ (at $m \approx n/\sqrt 2$). 
This means that $\approx 0.76555 \,n$ decimal, or $2.5431 \,n$ binary, leading significant digits 
have to cancel in the summation (\ref{ldd}) to yield the final $\L _n$.
Arbitrary-precision computing is thus mandatory.
This high instability in the specification of $\L_n$ (and already of $\l _n$, to a lesser degree) may be 
a price to pay for real-axis data that will signal phenomena located at very high imaginary parts.

\subsection{Analytical questions}

Fine-tuned as it is, the construction (\ref{ldd}) still has residual flexibility: 
e.g., in \cite[App.]{V2} we exhibit a more symmetrical - algebraically less elementary - variant 
(hence the plural in our main title). 
Inversely then, we may hope that simpler or better conditioned variants 
of (\ref{ldd}) could emerge in some future.

The right-hand side of (\ref{hf}) only shows the start of a double expansion 
in integer powers of $1/\log n$ and complex powers of $n$, thus constituting a \emph{transseries}
in the variable $1/\log n$ - whose analytical properties wholly remain to be investigated.

The asymptotic alternative (\ref{hf})--(\ref{ht}) entails an asymptotic Li-like criterion 
for the sequence $\{\L _n\}$: RH true $\iff$ for some $n_0$, $\L _n>0 \ (\forall n>n_0)$.
We cannot rule out $n_0=0$ (full Li's criterion), but we have no proof of this either.

\section{The Davenport--Heilbronn counterexamples}

They are ``twisted zeta functions" defined by the Dirichlet series \cite{DaH}\cite{BG}
\beq
\label{DHf}
f_\pm (x) \defi \frac{1}{1^x} + \frac{\tau _\pm}{2^x} - \frac{\tau _\pm}{3^x} - \frac{1}{4^x} + \frac{0}{5^x} + \cdots
\eeq
where $\tau _\pm = - \phi \pm \sqrt{1+\phi ^2}$ with $\phi = \hf (1+\sqrt 5)$ (the golden ratio),
and the numerators are repeated periodically mod 5.

For those specific $\tau_\pm$-values, $f_\pm (x)$ retain some properties like those of Riemann's $\z (x)$, 
e.g., countably many explicit values, and functional equations similar to (\ref{rfe}).
However, they lose the arithmetical properties of $\z (x)$ such as (\ref{epf}), 
and part of their (completed functions') zeros lie \emph{off} the critical line $L$.

\begin{figure}[t]
\center
\includegraphics[scale=.39]{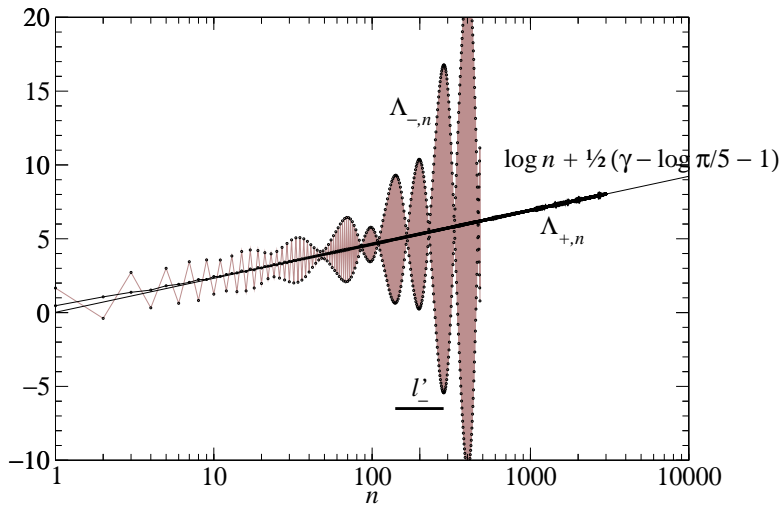}
\caption{\label{fg6} \small As Fig.~\ref{fg4}, but for the two sequences
$\{\L _{\pm,n}\}$ associated with the Davenport--Heilbronn functions (\ref{DHf}),
stopping $\L _{-,n}$ at $n=480$ before severe overflow occurs;
straight line: the generalized-RH-true asymptotic form (\ref{DHC}). 
The horizontal segment depicts the ($\log n$)-wavelength $l'_-$ 
of the leading oscillatory contribution $\propto n^{\rho '_- -1/2}$ 
to $\L _{-,n}$ (ignoring the factor $(-1)^n$ in (\ref{AS})): 
$l'_- \approx 2\pi /8.91836 \approx 0.70452$.}
\end{figure}

As with $\z (x)$, Keiper--Li sequences $\{\l _{\pm,n}\}$ can be defined for $f_\pm$,
and then, variants $\{\L _{\pm,n}\}$ in elementary closed form as well
(now using Bernoulli-polynomial values at 1/5 and 2/5).
These sequences $\{\L _{\pm,n}\}$ can then serve to numerically probe the RH-false branch 
(\ref{hf}) of our asymptotic alternative. \cite[\S~4.4]{V2}

The two functions $f_\pm$ yield contrasting numerical results (Fig.~\ref{fg6}).

\begin{figure}[h]
\center
\includegraphics[scale=.39]{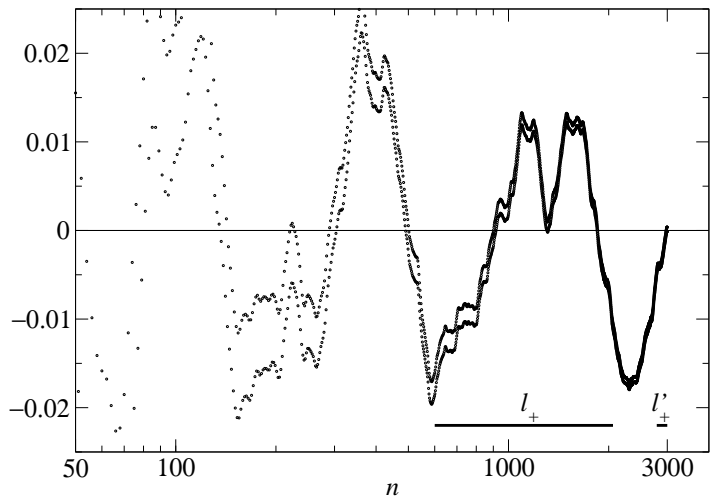}
\caption{\label{fg7} \small As Fig.~\ref{fg5}, but for the remainder $\delta \L _{+,n}$.
The horizontal segments depict the ($\log n$)-wavelengths $l_+ ,\ l'_+$ of the expressions $n^\rho $ 
for $\rho = \hf + 5.09416 \, \mi$ (the first zero of $f_+$, on $L$), and $\rho = \rho '_+$ 
(its first zero off $L$) respectively: $l_+ \approx 1.23341,\ l'_+ \approx 0.0733\,$.}
\end{figure}

For (the completed function of) $f _+$, its lowest-$T$ zero off the line $L$ is 
$\rho '_+ \approx 0.808517 + 85.699348 \, \mi$ \cite{S}.
To detect it through the sequence $\{ \L _{+,n} \}$, (\ref{TN}) gives a threshold 
$n \approx T^{1+2/t} \approx (85.7)^{7.48} \approx 3 \cdot 10^{14}$: 
accordingly, in our computed range $\L _{+,n}$ sticks to the RH-true prediction as generalized to this case,
\beq
\label{DHC}
\L _{\pm,n} \approx \log n + \hf (\g -\log \pi /5 -1) .
\eeq

Whereas for the case of $f _-$, the lowest-$T$ zero off $L$ is 
$\rho '_- \approx {2.30862 + 8.91836 \, \mi}$ 
(\cite{BS}, where our $f_-$ is denoted $f_2$).
To detect it through $\{ \L _{-,n} \}$, (\ref{TN}) now gives a threshold 
$n \approx T^{1+2/t} \approx (8.92)^{2.11} {\approx 100}$.
Indeed, $\L _{-,n}$ briefly starts along (\ref{DHC}) on average,
but the oscillating contribution like (\ref{AS}) from the zero $\rho '_-$ 
quickly turns dominant. This actually models what one should see
at much higher $n$ for $\{ \L _{+,n} \}$, 
and at some still higher $n$ for $\{ \L _n \}$ itself (the Riemann case) 
if RH is ultimately false.

On the other hand, $f_+$ is more suitable for practising to spot an early, hence weak, 
signal of RH-violation from $\rho '_+$ within the remainder $\delta \L _{+,n}$ in (\ref{DHC}).
(With $f_-$, the signal from $\rho '_-$ takes over too soon to allow that.) 
On $f_+$, the uncertainty-principle bound (\ref{ucp}) for the first zero $\rho '_+ \notin L $ 
gives $n \gtrsim 1100$: that leaves a large $n$-interval $\approx {[10^3,3 \cdot 10^{14}]}$ 
as training ground, to scan $(-1)^n \delta \L _{+,n}$ 
for an oscillation $\propto n^{\rho '_+ -1/2}/\log n$ 
(of ($\log n$)-wavelength: $l'_+ = 2\pi/\Im \rho '_+ \approx 0.0733$) (Fig.~\ref{fg7}). 
At our highest data point $n=4000$ its amplitude is $\approx 6 \cdot 10^{-5}$ by (\ref{AS}), still tiny, but it will grow like $n^{0.3085}/\log n$.
Then, the faster the background part of $\delta \L _{+,n}$ 
($=o(1)$, due to zeros of $f_+$ on $L$) would decrease,
the sooner that signal from $\rho '_+$ might stand out in the above interval. 
Such wishful thinking suffices to suggest what to seek next: 
higher-$n$ data for sure, 
but also stronger bounds on the remainder $\delta \L _{+,n}$ 
($\delta \L _n$ for the Riemann case, just as \cite{Lg} did with $\l _n$), 
and refined signal processings
- the end goal being to most efficiently use the sequence $\{\L _n\}$ itself for tests of RH.


\bigskip 

\begin{thebibliography}{99}

\bibitem{AR} J. Arias de Reyna, {\it Asymptotics of Keiper--Li coefficients\/},
Funct. Approx. Comment. Math. {\bf 45} (2011) 7--21.

\bibitem{A} A. Avila, {\it Convergence of an exact quantization scheme\/}, 
Commun. Math. Phys. {\bf 249} (2004) 305--318.

\bibitem{BS} E.P. Balanzario and J. S\'anchez-Ortiz,
{\it Zeros of the Davenport--Heilbronn counterexample\/},
Math. Comput. {\bf 76} (2007) 2045--2049.

\bibitem{BB} R. Balian and C. Bloch, {\it Solutions of the Schr\"odinger equation in terms of classical paths\/},
Ann. Phys. (NY) {\bf 85} (1974) 514--545.

\bibitem{BPV} R. Balian, G. Parisi, A. Voros, {\it Quartic oscillator\/}, 
in: S.~Albeverio et al. (eds.), {\it Feynman Path Integrals\/} [(Proceedings, Marseille 1978),
Lecture Notes in Physics {\bf 106}, Springer, Berlin (1979) 337--360
(errata: \cite[footnote p.~209]{Vz} or \cite[p.~203--204]{CC}).

\bibitem{BG} E. Bombieri and A. Ghosh,
{\it Around the Davenport--Heilbronn function\/}, 
Uspekhi Mat. Nauk {\bf 66} (2011) 15--66, Russian Math. Surveys {\bf 66} (2011) 221--270.

\bibitem{BL} E. Bombieri and J.C. Lagarias,
{\it Complements to Li's criterion for the Riemann Hypothesis\/},
J. Number Theory {\bf 77} (1999) 274--287.

\bibitem{CC} D. Chudnovsky and G. Chudnovsky (eds.), 
{\it The Riemann Problem, Complete Integrability and Arithmetic Applications\/}
(Proceedings, IH\'ES and Columbia University, 1979--1980),
Lecture Notes in Mathematics {\bf 925}, Springer, Berlin (1982).

\bibitem{C1} M.W. Coffey, 
{\it Toward verification of the Riemann Hypothesis:
application of the Li criterion\/}, Math. Phys. Anal. Geom. {\bf 8} (2005) 211--255.

\bibitem{C2} M.W. Coffey, 
{\it New results concerning power series expansions of the Riemann xi function and the Li/Keiper constants\/}, 
Proc. R. Soc. Lond. {\bf A 464} (2008) 711--731.

\bibitem{DaH} H. Davenport and H. Heilbronn, {\it On the zeros of certain Dirichlet series I, II\/},
J. London Math. Soc. {\bf 11} (1936) 181--185, 307--312.

\bibitem{Di} R.B. Dingle, 
{\it Asymptotic Expansions: their Derivation and Interpretation\/}, 
Academic Press, London (1973).

\bibitem{E} J. \'Ecalle, {\it Les fonctions r\'esurgentes I\/},
Publications Math\'ematiques d'Orsay {\bf 81-05} 
({\tt http://sites.mathdoc.fr/PMO/PDF/E{\detokenize{_}}ECALLE{\detokenize{_}}81{\detokenize{_}}05.pdf})
(and {\it II\/}, PMO {\bf 81-06} (1981); {\it III\/}, PMO {\bf 85-05} (1985)).

\bibitem{G} X. Gourdon, {\it The $10^{13}$ first zeros of the Riemann Zeta function,
and zeros computation at very large height\/}, preprint (Oct. 2004),\break
{\tt\footnotesize http://numbers.computation.free.fr/Constants/Miscellaneous/zetazeros1e13-1e24.pdf}

\bibitem{J} F. Johansson, {\it Rigorous high-precision computation of the Hurwitz zeta function
and its derivatives\/}, Numer. Algor. {\bf 69} (2015) 253--270.

\bibitem{KT} T. Kawai and Y. Takei, {\it Algebraic Analysis of Singular Perturbation Theory\/},
Translations of Mathematical Monographs {\bf 227}, Amer. Math. Soc. (2005) 
[Japanese: Iwanami shoten (1998)], and refs. therein.

\bibitem{K} J.B. Keiper, {\it Power series expansions of Riemann's $\xi$ function\/},
Math. Comput. {\bf 58} (1992) 765--773.

\bibitem{KS} T. Koike and R. Sch\"afke, in preparation.

\bibitem{Lg} J.C. Lagarias, {\it Li coefficients for automorphic $L$-functions\/}, 
Ann. Inst. Fourier, Grenoble {\bf 57} (2007) 1689--1740.

\bibitem{LI1} X.-J. Li,
{\it The positivity of a sequence of numbers and the Riemann Hypothesis\/},
J. Number Theory {\bf 65} (1997) 325--333.

\bibitem{M1} K. Ma\'slanka, {\it Li's criterion for the Riemann hypothesis --
numerical approach\/}, Opuscula Math. {\bf 24} (2004) 103--114.

\bibitem{MP} S. Minakshisundaram, \r{A}. Pleijel, 
{\it Some properties of the eigenfunctions of the Laplace-operator on Riemannian manifolds\/}, 
Can. J. Math. {\bf 1} (1949) 242--256.

\bibitem{Mi} G. Misguich, calculations for $n>20000$, using {\tt http://www.mpfr.org/} 
(private communications, 2017).

\bibitem{O} J. Oesterl\'e, {\it R\'egions sans z\'eros de la fonction z\^eta de Riemann\/},
typescript (2000, revised 2001, uncirculated).

\bibitem{P} G. Parisi, {\it Trace identities for the Schr\"odinger operator and the WKB method\/},
Preprint LPTENS 78/9 (\'Ecole Normale Sup\'erieure, Paris, March 1978), in \cite[p.178--183]{CC}.

\bibitem{Ri} B. Riemann, {\it \"Uber die Anzahl der Primzahlen unter einer gegebenen Gr\"osse\/},
Monatsb. Preuss. Akad. Wiss. (Nov. 1859) 671--680;
English translation, by R.~Baker, Ch.~Christenson and H.~Orde: 
{\it Bernhard Riemann: Collected Papers\/}, paper VII, Kendrick Press, Heber City, UT (2004) 135--143.

\bibitem{SKK} M. Sato, T. Kawai and M. Kashiwara,
{\it Microfunctions and pseudo-differential equations\/}, in: H.~Komatsu (ed.),
{\it Hyperfunctions and pseudo-differential equations\/} (Proceedings, Katata 1971),
Lecture Notes in Mathematics {\bf 287}, Springer, Berlin (1973) 265--529.

\bibitem{Se} S.K. Sekatskii, {\it Generalized Bombieri--Lagarias' theorem and 
generalized Li's criterion with its arithmetic interpretation\/},
Ukr. Mat. Zh. {\bf 66} (2014) 371--383, Ukr. Math. J. {\bf 66} (2014) 415--431.

\bibitem{S} R. Spira, {\it Some zeros of the Titchmarsh counterexample\/}, Math. Comput. {\bf 63} (1994) 747--748.

\bibitem{Ti} E.C. Titchmarsh, {\it The Theory of the Riemann Zeta-Function\/},
2nd ed. revised by D.R. Heath-Brown, Oxford Univ. Press, Oxford (1986).

\bibitem{VGS} A. Voros, {\it Oscillateur quartique et m\'ethodes semi-classiques\/},
in: {\it S\'eminaire Goulaouic--Schwartz 1979--1980\/}, \S~VI (nov. 1979)
({\tt https://eudml.org/doc/111764} with errata in \S~``Notes").

\bibitem{Vz} A. Voros, {\it The zeta function of the quartic oscillator\/}, 
Nucl. Phys. {\bf B165} (1980) 209--236 (errata :
{\tt www.ipht.fr/Docspht//articles/t79/046/public/erratum.pdf}), augmented version: {\it Zeta functions of the quartic 
(and homogeneous anharmonic) oscillators\/}, in \cite{CC} p.184--208.

\bibitem{Vc} A. Voros, {\it Spectre de l'\'equation de Schr\"odinger et m\'ethode BKW\/},
Publications Math\'ematiques d'Orsay {\bf 81-09} (1981) (corrected copy:
{\tt https://www.ipht.fr/Docspht//articles/t81/120/public/81{\detokenize{_}}120.pdf}).

\bibitem{V0} A. Voros, {\it Correspondance semi-classique et r\'esultats exacts :
cas des spectres d'op\'erateurs de Schr\"odinger homog\`enes\/},
C.R. Acad. Sci., Paris, S\'er.~I {\bf 293} (1981) 709--712
({\tt https://gallica.bnf.fr/ark:/12148/bpt6k64461526/f351.image}),
English translation: {\it Semiclassical correspondence and exact results:
the case of the spectra of homogeneous Schr\"odinger operators\/},
J. Physique Lettres {\bf 43} (1982) L-1--L-4 (erratum: ibid. p.159).

\bibitem{VB} A. Voros, {\it Le probl\`eme spectral de Sturm--Liouville : le cas de l'oscillateur quartique\/},
(expos\'e no.602, Nov. 1982), in: S\'eminaire Bourbaki {\bf 25} (1982--1983), 
Ast\'erisque {\bf 105--106}, Soc. Math. France (1983) 95--104 
({\tt eudml.org/doc/110019} with reference updates in \S~``Notes").

\bibitem{V1} A. Voros, {\it The return of the quartic oscillator. The complex WKB method\/},
Ann. Inst. H. Poincar\'e {\bf A39} (1983) 211--338 ({\tt eudml.org/doc/76217} with errata in \S~``Notes").

\bibitem{Ve} A. Voros, {\it The general 1D Schr\"odinger equation as an exactly solvable problem\/}, in: 
Y.~Takei (ed.), {\it Recent trends in exponential asymptotics\/} (Proceedings, Kyoto 2004),
RIMS K\^oky\^uroku {\bf 1424} (2005) 214--231.

\bibitem{VK} A. Voros, {\it From exact-WKB toward singular quantum perturbation theory~II\/}, 
in: T.~Aoki {\it et al.\/} (eds.), {\it Algebraic analysis of differential equations\/}
(Festschrift in honor of T.~Kawai, Kyoto 2005), Springer, Tokyo (2008) 321--334.

\bibitem{V} A. Voros, {\it Sharpenings of Li's criterion for the Riemann Hypothesis\/},
Math. Phys. Anal. Geom. {\bf 9} (2006) 53--63
(erratum: our asymptotic forms for $\l _n$ in the RH false case have wrong sign).

\bibitem{V2} A. Voros, 
{\it Discretized Keiper/Li approach to the Riemann Hypothesis\/},
Exp. Math. {\bf 29} (2020) 452--469 (publ. online Jul.~17, 2018).

\bibitem{W} S. Wolfram, Mathematica, 3rd ed., Wolfram Media/Cambridge University Press, New York (1996).

\end{thebibliography}
\end{document}